\documentclass[12pt]{amsart}
\usepackage{amssymb}
\usepackage[dvips]{graphics}
\ifx\mathrm\undefined\let\mathrm\rm\fi
\ifx\mathbf\undefined\let\mathbf\bf\fi
\ifx\mathfrak\undefined\let\mathfrak\frak\fi
\ifx\mathcal\undefined\let\mathcal\cal\fi
\ifx\mathbb\undefined\let\mathbb\Bbb\fi
\ifx\emph\undefined\let\emph\it\fi
 at9.98pt
\newcommand{\Z}{{\mathbb Z}}
\newcommand{\ZZ}{{\mathbb Z_{>0}}} 

\newcommand{\R}{{\mathbb R}}
\newcommand{\C}{{\mathbb C}}

\newcommand{\A}{\mathcal{A}}
\newcommand{\PP}{\mathcal{P}}  
\newcommand{\dontprint}[1]
{\relax}
\newtheorem%
{thm}{Theorem}
\newtheorem%
{proposition}
{Proposition}
\newtheorem%
{lemma}
{Lemma}
\newtheorem%
{lemmadef}[thm]{Lemma-Definition}
\newtheorem%
{corollary}
{Corollary}
\newtheorem%     
{conjecture}
{Conjecture}

\begin{document}
\title{ Critical points of functions, $sl_2$ representations,\\
 and Fuchsian differential equations\\ with only univalued solutions}
%\def\shorttitle{ Critical points, $sl_2$ representations, 
%and Fuchsian equations}
%*********************
\author[{}]
{I. Scherbak${}^{\ast}$
\and A. Varchenko${}^{\dagger,1}$}
\thanks{${}^1$Supported in part by NSF grant  DMS-9801582}
\maketitle
\medskip 
\centerline{\it ${}^{\ast}$School of Mathematical Sciences,
Tel Aviv University,}
\centerline{\it Ramat Aviv, Tel Aviv 69978, Israel} 
\medskip
\centerline{\it ${}^{\dagger}$Department of Mathematics, University of
  North Carolina at Chapel Hill,} 
\centerline{\it Chapel Hill, NC 27599-3250, USA} 
\medskip
 
%********************
\pagestyle{myheadings}
\markboth{I. Scherbak and A. Varchenko}
{ Critical points, $sl_2$ representations,
and Fuchsian equations}       
%*********************                                        
%\tableofcontents

\bigskip
\centerline{\emph{Dedicated to V.~I.~Arnold on his 65th birthday}}

\begin{abstract}
${}$
Let a second order Fuchsian differential equation with only univalued
solutions have finite singular points at $z_1,...,z_n$ with exponents
$(\rho_{1,1},\rho_{2,1})$,...,$(\rho_{1,n},\rho_{2,n})$.
Let the exponents at infinity be $(\rho_{1,\infty},\rho_{2,\infty})$.
%Let the point at infinity have
%exponents  $(\rho_{1,\infty},\rho_{2,\infty})$.
Then for fixed generic $z_1,...,z_n$,
the number of such Fuchsian equations is equal to the multiplicity of the
irreducible $sl_2$ representation of dimension
$| \rho_{2,\infty} - \rho_{1,\infty} |$
in the tensor product of irreducible $sl_2$ representations of dimensions
$| \rho_{2,1} - \rho_{1,1} |$, ... ,
$| \rho_{2,n} - \rho_{1,n} |$.  
 To show this we count the number of critical
points 
of a suitable function which plays the crucial role in constructions
of the hypergeometric solutions of the $sl_2$ KZ equation and of the Bethe 
vectors in the $sl_2$ Gaudin model. As a byproduct of this study we
conclude 
that the set of Bethe vectors is a basis in the space of states for the $sl_2$ 
inhomogeneous Gaudin model.

\end{abstract}

\section {Introduction}

\subsection{Critical points and $sl_2$ representations}
 
Consider the Lie algebra $sl_2$ with standard generators
$e\,,f\,,h\,$,  $[e,f]=h,\ \ [h,e]=2e,\ \ [h,f]=-2f\,$.
Let $L_a$ be the irreducible $sl_2$ module with highest
weight $a\in \mathbb C\,$. The module $L_a$ is generated by
its singular vector $v_a\,, ev_a=0\,, hv_a=av_a\,$.
Vectors $v_a\,, fv_a\,, f^2v_a\,,\dots $
form a basis of $L_a\,$. If  $a$ is a nonnegative integer,
then $\dim L_a=a+1$; otherwise $L_a$ is infinite-dimensional.
 
If $m_1,\dots,m_n\,$ are nonnegative integers, then the tensor product
$L^{\otimes m}=L_{m_1}\otimes\cdots \otimes L_{m_n}$
is a direct sum of irreducible representations with highest
weights $l(m)-2k\,$, where 
$$
l(m)\ =\ m_1 \ +\ ... \ + \ m_n
$$ 
and $k$
is a nonnegative integer.
Let $w(m,k)$ be the multiplicity of $L_{l(m)-2k}$ in $L^{\otimes m}\,$.
We have
$$
w(m,k)\geq 0\ \ {\rm if}\ \  l(m)-2k\geq 0 \,;\quad
w(m,k)=0\ \ {\rm if}\ \ l(m)-2k<0 \,.
$$
  
Let $z=(z_1,\dots, z_n) \in \mathbb C^n$ be a point
with pairwise distinct coordinates. Let
$$
{\mathcal A}={\mathcal A}_{k,n}(z)=\bigcup_{i=1}^k
\bigcup_{l=1}^n\{t\in\mathbb C^k\ \vert\ t_i=z_l\}
\bigcup_{1\leq i<j\leq k}\{t\in\mathbb C^k\ \vert\ t_i=t_j\}
$$
be a (discriminantal) arrangement of hyperplanes in $\mathbb C^k\,$,
and ${\mathcal C}={\mathcal C}_{k,n}(z)\,$ its complement. 
For $m=(m_1,\dots,m_n)\in \mathbb C^n\,$, consider the
multivalued function $\Phi:\ {\mathcal C}\rightarrow\mathbb C\,$,
$$
\Phi_{k,n}(t)=\Phi_{k,n}(t;z, m)=\prod_{i=1}^k\prod_{l=1}^n
(t_i-z_l)^{-m_l}\prod_{1\leq i<j\leq k}(t_i-t_j)^2 \,.
$$
A point $t^0\in {\mathcal C}$ is
called {\it a critical point of} $\Phi$ if
$$
\frac{\partial\Phi}{\partial t_i}(t^0)=0,\quad
i=1,\dots, k \,.
$$
The symmetric group $S^k$ acts on ${\mathcal C}$ permuting coordinates.
Each orbit consists of $k!$ points. The action preserves the
critical set of the function $\Phi_{k,n}(t)\,$.
 
Let $\lambda_1=\sum t_i,\ \ \lambda_2=\sum t_it_j,\ \ \dots,
\ \ \lambda_k=t_1 \cdots  t_k$ be the standard symmetric
functions of $t_1, ... , t_k\,$.
Denote $\C^k_\lambda$ the space with coordinates $\lambda_1,...,\lambda_k$.
 
Our first main result is 
 
\begin{thm}\label{TM} \ \ 
Let $m_1,\dots,m_n\, \in \Z_{\geq 0}$  and $k \in \Z_{>0}$.
\begin{itemize}
\item
If $l(m)+1- k> k$, then for generic $z$ all critical points of
$\Phi_{k,n}(t)$ are nondegenerate and the critical set consists
of $w(m,k)$ orbits.
\item
If $l(m)+1-k = k$, then for any $z$ the function
$\Phi_{k,n}(t)$ does not have critical points.
\item
If $\ 0 \leq l(m)+1- k < k$,  then for generic $z$
the function $\Phi_{k,n}(t)$ can have  only non-isolated critical points.
Written in symmetric coordinates $\lambda_1,...,\lambda_k\,$,
the critical set consists of $w(m,l(m)+1-k)$ straight lines
in the space $\C^k_\lambda\,$.
\item
If $l(m)+1-k < 0$, then for any $z$ the function
$\Phi_{k,n}(t)$ does not have critical points.
\end{itemize}
\end{thm}
 
In this paper, the words ``a point $z $ is generic''
mean that $z$ does not belong to a suitable proper algebraic subset in $\C^n$.

\medskip\noindent
{\bf Remark.} Assume that $m_1, ... \,, \,m_p \in \Z_{>0}$ and $
m_{p+1} = ... = m_n = 0$. Then for any $k$, we have
$\Phi_{k,n}(t_1, ... , t_k; z_1, ... , z_n,\, m_1, ... , m_n)
=
\Phi_{k,p}(t_1, ... , t_k; z_1, ... , z_p,\, m_1, ... , m_p)$.
For generic $z_{p+1}, ... , z_n$, the two functions have the same number of isolated 
critical points and the same number of critical curves. We also have
$w( \,(m_1, ... , m_n),\, k\, ) = w( \,(m_1, ... , m_p),\, k\,)$.
Thus to prove the Theorem it is enough to consider the case when all $m_1, ... , m_n$
are positive integers.

\medskip\noindent
{\bf Example.}\ \  Let $n=2$ and  $z=(0,1)\,$. We have
$$
\Phi_{k,2}(t)=\prod_{i=1}^k
t_i^{-m_1}(t_i-1)^{-m_2}\prod_{1\leq i<j\leq k}(t_i-t_j)^2 \,.
$$
The critical point system of the function
$\Phi_{k,2}(t)$ being written
with respect to symmetric coordinates
$\lambda_1,\dots,\lambda_k$ is the following linear system,
\begin{eqnarray*}
(p+1)(-m_1+p)\lambda_{k-p-1}=(k-p)(-m_1-m_2+k+p-1)\lambda_{k-p} \,,
\end{eqnarray*}
where $ p=0,\dots,k-1\,$, and $\lambda_0=1\,$, see Lemma 1.3.4 in \cite{V}.

For $m_1,m_2\in\Z_{>0}\,$, there are four possibilities.    

\medskip\noindent
{\it (i)}\ \ If $k\leq m_1,m_2\,$, then the linear system has
a single solution which defines $k!$ nondegenerate critical points
of $\Phi_{k,2}(t)\,$. In this case the multiplicity of $L_{m_1+m_2-2k}$
in $L_{m_1}\otimes L_{m_2}$ is $w(m,k)=1\,$.
 
\medskip\noindent
{\it (ii)}\ \ If $k$ is greater than exactly one of the numbers
$m_1,m_2\,$, then the linear system  still has a single solution,
but the solution defines points lying in the arrangement ${\A}\,$.
This means that $\Phi_{k,2}(t)$ does
not have critical points. In this case $w(m,k) = w(m, l(m)+1-k)=0$.
 
\medskip\noindent
{\it (iii)}\ \ If $m_1,m_2< k \leq m_1+m_2 +1\,$, then the rank of
the linear system  is $k-1\,$. The solutions form a straight line
in the space $\C^k_\lambda$ with coordinates
$\lambda_1,\dots,\lambda_k\,$. The line defines a curve
of critical points of the function $\Phi_{k,2}(t)\,$.
In this case  $w(m,l(m)+1-k)=1\,$.
 
\medskip\noindent
{\it (iv)}\ \ If $m_1+m_2+1< k\,$, then the system  again has a
single solution which defines points lying in the arrangement
$\A\,$. The function $\Phi_{k,2}(t)$ does
not have critical points. 

\bigskip
For negative exponents $m_1,...,m_n$ and real $z_1, ... , z_n$,
the function $\Phi_{k,n}(t;z, m)$
has only nondegenerate critical points and the critical set consists of
${k+n-2\choose n-2}$ orbits \cite{V}. 
If the exponents tend to positive integer values
so that $l(m)-2k$ remains nonnegative, 
some of critical points vanish at edges of the arrangement 
${\mathcal A}_{k,n}(z)$. 
To prove the first part of Theorem \ref{TM} we count the number of 
vanishing critical points.

To show that critical points of $\Phi_{k,n}(t;z, m)$ form lines, 
if $l(m)-2k$ is negative,
we use the connection of critical points with Fuchsian equations having
polynomial solutions.

On the number of critical points of a product of generic powers of 
arbitrary linear functions see  \cite{V, OT, Si}. In that case of generic exponents
the critical points are isolated and nondegenerate and their number is equal to the absolute
value of the Euler characteristic of the complement to the arrangement  of zero hyperplanes of
the linear functions. In contrast to generic exponents, the exponents of the function
$\Phi_{k,n}(t;z, m)$ in  Theorem \ref{TM}
are highly resonant. It would be very interesting to find out how much
of the phenomenon described in Theorem \ref{TM} can be generalized to
more general arrangements.

\subsection{Critical points and Fuchsian equations 
with polynomial solutions}

Consider the differential equation
\begin{eqnarray}\label{e}
u''(x)+p(x)u'(x)+q(x)u(x)=0\,
\end{eqnarray}
with meromorphic $p(x)$ and $q(x)$. A point $z_0\in\C$ is 
{\it an ordinary point} of the equation if
the functions $p(x)$ and $q(x)$ are holomorphic at $x=z_0\,$.
A non-ordinary point  is called {\it singular}.
 
The point $z_0\in\C$ is {\it a regular singular point} of the equation if
$z_0$ is a singular point,
$p(x)$ has a  pole at $z_0$ of order not greater than 1, and
$q(x)$ has a  pole at $z_0$ of order not greater than 2.

The equation has {\it an ordinary (resp., regular singular) 
point at infinity} if after the change
$x=1/\xi$ the point $\xi=0$ is an ordinary (resp., regular singular)
point of the transformed equation. 
%The equation   has a {\it regular singular point at infinity}
%if after the change $x=1/\xi$ the point $\xi=0$ is a regular
%singular point of the transformed equation.

Let $x=z_0$ be a regular singular point,
$$
p(x)=\sum_{l=0}^{\infty}{p_l}(x-z_0)^{l-1}\,,\quad
q(x)=\sum_{l=0}^{\infty}{q_l}(x-z_0)^{l-2}\,
$$
the Laurent series at $z_0\,$.
If the function
\begin{eqnarray}\label{X}
u(x)=(x-z_0)^{\rho}\sum_{l=0}^{\infty}c_l(x-z_0)^l,\quad c_0=1,
\end{eqnarray}
is a solution of  equation (\ref{e}), then $\rho$ must be 
a root of {\it the indicial equation}
$$
\rho^2+(p_0-1)\rho +q_0=0\,.
$$
The roots of the indicial equation are called
{\it the exponents} of the equation  at $z_0\,$.
 
If the difference $\rho_1-\rho_2$ of roots is not an integer, then
the equation has solutions of the form (\ref{X}) with $\rho=\rho_j\,$,
$j=1, 2\,$. If the difference $\rho_1-\rho_2$ is a nonnegative integer, 
then the equation has a solution $u_1$ of the form (\ref{X}) 
with $\rho=\rho_1\,$.
The second linearly independent solution $u_2$ is either
of the form
$$
u_2(x)=(x-z_0)^{\rho_2}\sum_{l=0}^{\infty}d_l(x-z_0)^l\,,\quad d_0=1\,,
$$
or
$$
u_2(x)=u_1(x)\ln (x-z_0) + 
(x-z_0)^{\rho_2}\sum_{l=0}^{\infty}d_l(x-z_0)^l\,.
$$

A differential equation  with only regular singular points 
is called {\it Fuchsian}. Let the singular points of a Fuchsian 
equation be $z_1,\dots,z_n$ and infinity.
Let $\rho_{1,j}$ and $\rho_{2,j}$ be the exponents at 
$z_j,\ 1\leq j\leq n,$ and $\rho_{1,\infty},\ \rho_{2,\infty}$ 
the exponents at infinity. Then
$$
\rho_{1,\infty}+\rho_{2,\infty}+\sum_{j=1}^n(\rho_{1,j}+\rho_{2,j})=n-1\,.
$$
 
Consider the equation
\begin{eqnarray}\label{Ee}
F(x)u''(x)+G(x)u'(x)+H(x)u(x)=0\,,
\end{eqnarray}
where $F(x)$ is a polynomial of degree $n$, and $G(x)\,,\ H(x)$ 
are polynomials  of degree not greater than   $n-1\,,\ n-2\,,$ respectively.
If $F(x)$ has no multiple roots, then the equation is Fuchsian.
Write
\begin{eqnarray}\label{M} 
F(x) = \prod_{j=1}^n (x-z_j)\,,\qquad
\frac{G(x)}{F(x)} = \sum_{j=1}^n\frac{-m_j}{x-z_j}\, 
\end{eqnarray}    
for suitable complex numbers $m_j\,, z_j\,$. Then $0\,$ and $m_j+1\,$ 
are exponents at $z_j$ of  equation (\ref{Ee}). 
If $-k$ is one of the exponents at $\infty\,$, then 
the other is $k-l(m)-1$.

\medskip\noindent
{\bf Problem} \cite{S}, Ch.~6.8.\ \ 
{\it Given polynomials $F(x)\,,\ G(x)$ as above,
\begin{itemize}
\item[(i)] find a polynomial $H(x)$
of degree at most $n-2$ such that  equation (\ref{Ee})
has a polynomial solution of a preassigned degree $k\,$;
\item[(ii)] 
find the number of solutions to   Problem (i).
\end{itemize}} 

The following result is classical.

\begin{thm}\label{PS}{\rm Cf. \cite{S}, Ch.~6.8.}
\begin{itemize}
\item
Let $u(x)$ be a polynomial
solution of (\ref{Ee}) of degree $k$ with roots $t_1^0,..., t_k^0\,$
of multiplicity one. Then $t^0=(t_1^0,\dots,t_k^0)\,$ is a critical 
point of the function
$\Phi_{k,n}(t;z,m)\,$, where $z=(z_1,\dots,z_n)$ and $m=(m_1,\dots,m_n)\,$.
\item Let $t^0$ be a critical point of the function $\Phi_{k,n}(t;z,m)\,$, 
then the polynomial $u(x)$ of degree $k$ with
roots $t_1^0,..., t_k^0$ is a solution of (\ref{Ee}) with
$H(x)=(-F(x)u''(x)-G(x)u'(x))/u(x)$
being  a polynomial of degree at most $n-2\,$.
\end{itemize}
\end{thm}

A critical point of the function $\Phi_{k,n}(t;z,m)$ defines a Fuchsian
   differential equation and its polynomial solution.
   The Fuchsian differential equation defined by a critical point $t^0$
   will be called {\it associated} and denoted $E(t^0,z,m)$.

According to Theorem \ref{PS}, the orbits of critical points of $\Phi_{k,n}(t;z,m)$ label solutions
   to Problem (i), and Problem (ii) turns out to be the question on the number
   of the orbits of critical points of $\Phi_{k,n}(t;z,m)$.

   For fixed real $z_1,\dots, z_n$ and negative $m_1,\dots, m_n$,
   Problem (ii) was solved in the 19th century by Heine and Stieltjes.
They showed that under these conditions the number of solutions is equal to
${k+n-2\choose n-2}$, see \cite{S},  Ch.~6.8.

\subsection{Fuchsian equations with only polynomial solutions}\label{diff-eqn}

If all solutions of the Fuchsian equation (\ref{Ee}) are polynomials,
then the numbers $m_1,...,m_n$ in (\ref{M}) are nonnegative integers. 
If $k$ is the degree of the generic polynomial
solution of that equation, then $k> l(m)+1-k$ and the equation 
also has polynomial solutions of degree $l(m)+1-k$.

Assume that all solutions of equation  (\ref{Ee}) are polynomials.
If $m_j=0$ for some $j$, then $x=z_j$ is a regular point of the equation.
Indeed, the function $G(x)$ is clearly divisible by $x-z_j$. We also have
$$
{H(x) \over F(x)}\ = \ 
- \ {u''(x) \over u(x) }\ -\
{G(x) \over F(x)} \ 
{u'(x) \over u(x) }\ 
$$
for any solution $u(x)$. Hence $H(x)$ is divisible by $x-z_j$.

Our second main result counts for
generic $z_1,...,z_n$ the number of Fuchsian equations with
fixed positive integers $m_1,...,m_n, \, k$ having only polynomial solutions.

\begin{thm}\label{main2}
Let $z_1,\dots,z_n$ be pairwise distinct complex numbers, let 
$m_1,\dots,m_n, \, k\in \Z_{>0}$.
Let $F(x)$ and  $G(x)$ be determined by (\ref{M}).
\begin{itemize}
\item
If $k > l(m)+1-k \geq 0$,
then for generic $z_1,...,z_n$ there exist exactly 
$w(m,l(m)+1-k)$ polynomials
$H(x)$ of degree not greater than $n-2$ such that all solutions of
 equation (\ref{Ee})
are polynomials with the degree of the generic solution 
equal to $k$.
\item
If $k \leq l(m)+1-k$ or $l(m)+1-k < 0$ 
then for any $z_1,...,z_n$ there are no such polynomials $H(x)$.
\end{itemize} 
\end{thm}

The assumption that $z_1, ... ,z_n$ are generic is essential.
 
\medskip\noindent
{\bf Example.} \ \ Let $m=(1,1,1,)$ and $k=3$. Then $l(m)+1-k=1$ and
$w(m,l(m)+1-k)=2$. Let $z=(0,1,c)$,
$$
F(x)=x(x-1)(x-c)\,,\quad
G(x)=\left(-\frac1x-\frac1{x-1}-\frac1{x-c}\right)\cdot F(x)\,.
$$
Let $v(x)=x-\alpha$ be a solution of equation (\ref{Ee})
with these $F(x)$ and $G(x)\,$ and a suitable $H(x)$. 
The number $\alpha$
is a critical point of the function
$\Phi(x)=x^{-1}(x-1)^{-1}(x-c)^{-1}\,$, i.e. a root of the quadratic
equation $3x^2-2(c+1)x+c=0\,$.
 
If the discriminant of this equation
$\Delta=c^2-c+1$ is non-zero, then there are two distinct roots and
there are two linear polynomials
$$
H(x)=-\frac{G(x)}{x-\alpha}
$$
such that the corresponding equation (\ref{Ee}) has only polynomial solutions
with cubic generic solutions.
 
If the discriminant vanishes, i.e. the numbers $0, 1, c$ form an equilateral
triangle, then there is only one critical point and only one differential equation.

\medskip\noindent
{\bf Corollary of Theorem \ref{main2}. } 

{\it
Let all solutions of a second order Fuchsian differential equation
   be univalued. Let the singular points be $z_1,...,z_n$  and infinity.
   Let $\rho_{1,j}\,$ and $\rho_{2,j}\,$ be the exponents at $z_j$, and
   $\rho_{1,\infty}\,,\,\rho_{2,\infty}$ the exponents at infinity. 
Then for fixed generic $z_1,...,z_n$,
the number of such Fuchsian equations is equal to the multiplicity of the
irreducible $sl_2$ representation of dimension
$| \rho_{2,\infty} - \rho_{1,\infty} |$
in the tensor product of irreducible $sl_2$ representations of dimensions
$| \rho_{2,1} - \rho_{1,1} |$, ... ,
$| \rho_{2,n} - \rho_{1,n} |$.  }

\medskip\noindent
{\bf Proof.}\ \
Let
$$
v''(x)\ +\ p(x) v'(x)\ +\ q(x) v(x)\ =\ 0\,
$$ 
be such an equation. Assume that 
$\rho_{1,1}\,<\,\rho_{2,1}$, ..., $\rho_{1,n}\,<\,\rho_{2,n}$,
 $\rho_{1,\infty } \,<\, \rho_{2,\infty}\,$.
Change the variable, \ $v(x) = u(x) \ (x-z_1)^{\rho_{1,1}} \, \cdots \,
(x-z_n)^{\rho_{1,n}}$.  The new equation with respect to $u(x)$ is a Fuchsian 
differential equation of type ( \ref{Ee} ) with only polynomial solutions. Its exponents
are
$(0,\,\rho_{2,1} - \rho_{1,1})$, \ ... \ , \ $(0,\, \rho_{2,n} - \rho_{1,n})$, \
$( \rho_{1, \infty} + \sum_{j=1}^n \rho_{1,j}, \
\rho_{2, \infty} + \sum_{j=1}^n \rho_{1,j})$. By Theorem \ref{main2} for fixed generic
$z_1, ... , z_n, $ the number of such equations is $w(m, l(m)+1-k)$ where
$m = ( \rho_{2,1} - \rho_{1,1} - 1, \ ... \ ,\, \rho_{2,n} - \rho_{1,n} - 1 )$ and
$k = - \rho_{1, \infty} - \sum_{j=1}^n \rho_{1,j} $.
This gives the statement of the Corollary.
\hfill   $\triangleleft$  

\subsection{Two-dimensional spaces of polynomials with prescribed singularities}

Theorems \ref{TM} and \ref{main2} give the following corollary which can be considered 
as a statement from enumerative algebraic geometry, see for instance  \cite{GH}. Namely, 
consider a two-dimensional space $V$ of polynomials of one variable $x$ with complex coefficients.
Let $k_1$ be the degree of generic polynomials in $V$ and $k_2$ the degree of special polynomials
in $V$, $k_1 > k_2$.

For two functions $f(x), g(x)$ let $W(f,g)(x) = f'(x)g(x) - f(x)g'(x)$ be the Wronskian.
If  $f, g$ is a basis in $V$, then the Wronskian  has degree $k_1 + k_2 -1$ and 
does not depend on the choice of the basis up to multiplication by a nonzero constant.
The corresponding monic polynomial will be called the Wronskian of the space  and denoted
$W_V(x)$. Let
$$
W_V(x)\ =\ \prod_{l=1}^n\ (x - z_l)^{m_l}\ .
$$

We say that the vector space $V$ is nondegenerate if for any complex number $x_0$ there is a polynomial
$f(x)$ in $V$, such that $f(x_0)$ is not zero, and if the set of roots of a polynomial in 
$V$ of degree $k_2$ does not intersect the set  $z_1, ... , z_n$.

\medskip\noindent
{\bf Problem} 
{\it Assume that $k_1 > k_2$ and 
$W_V(x)\ =\ \prod_{l=1}^n\ (x - z_l)^{m_l}$ are fixed, $m_1 + ... + m_n$ 
\newline $ = k_1 + k_2 - 1$. 
What is the number of nondegenerate
vector spaces $V$ with such characteristics ? }

\medskip\noindent
{\bf Corollary of Theorems \ref{TM} and \ref{main2}. } 

{\it For  generic $z_1, ... , z_n$ the number of nondegenerate vector spaces
$V$ with such data is equal to the multiplicity of the representation
$L_{k_1-k_2-1}$ in the tensor product $L_{m_1} \otimes ... \otimes L_{m_n}$ .}

\subsection{Critical points and Bethe vectors }

For a positive integer $a$, let $L_a$ be the irreducible $sl_2$ module with 
highest weight $a$. 
The Shapovalov form on $L_a$ is the unique  
symmetric bilinear form $S_a$  such that
$S_a(v_a, v_a)=1$ and $S_a( ex, y) = S_a(x, fy)$ for all $x, y \in L_a$.

Let $\Omega={1\over 2} h\otimes h + e\otimes f + f \otimes e$ be the 
Casimir operator. 

For positive integers $m_1,...,m_n$, 
define on
$L^{\otimes m}=L_{m_1}\otimes\cdots \otimes L_{m_n}$
the Shapovalov form as $S=S_{m_1}\otimes\cdots \otimes S_{m_n}$.

For pairwise distinct complex numbers $z_1,...,z_n$ 
and any $i=1,...,n$, introduce a linear operator 
$H_i(z) : L^{\otimes m} \to L^{\otimes m}$,
\begin{eqnarray}\label{H} 
H_i(z)\,=\, \sum_{j,\, j\neq i} \,{\Omega^{(i,j)} \over z_i - z_j}\,.
\end{eqnarray}
Here  $\Omega^{(i,j)}$ is the operator acting as $\Omega$ in the $i$-th and
$j$-th factors and as the identity in all other factors of the tensor product. 
The operators $H_i(z), \, i=1,...,n,$ commute and are called the 
Hamiltonians of the Gaudin model
of an inhomogeneous magnetic chain \cite{G}. 

The  Bethe ansatz is a certain construction of eigenvectors for 
a system of commuting operators. The idea of the construction is to find a 
vector-valued function of a special form and determine its arguments 
in such a way that the value of this function is an eigenvector. 
The equations which determine the special values of arguments are called
the Bethe equations. The corresponding eigenvectors are called 
the Bethe vectors. The main problem of the Bethe ansatz is to show 
that the construction gives a basis of eigenvectors.
On the Bethe ansatz see for instance \cite{F, TV}.

One of the systems of commuting operators diagonalized by 
the Bethe ansatz is the system of Hamiltonians of the
Gaudin model, see \cite{G}.

Let $k$ be a positive integer. 
Let $J=(j_1,...,j_n)$ be a vector with  integer coordinates
such that $j_1+...+j_n=k$ and for any $l$ we have $0\leq j_l \leq m_l$. 
Introduce a vector in the tensor product, 
$f_Jv=f^{j_1}v_{m_1}\otimes ... \otimes f^{j_n}v_{m_n}$. 
Introduce a function 
$$
A_J(t_1,...,t_k, z_1,...,z_n)\,=\,
\sum_{\sigma\in \Sigma(k;\, j_1,...,j_n)} {}\,
\prod_{i=1}^k\, {1\over t_i - z_{\sigma(i)}}\,,
$$
the sum is over the set $\Sigma(k;\, j_1,...,j_n)$ of maps 
$\sigma$ from $\{1,...,k\}$ to $\{1,...,n\}$ such 
that for every $l$ the cardinality of $\sigma^{-1}(l)$ is equal to  $j_l$.

\begin{thm}\label{ReV}{\rm  \cite{RV, V}.}
\begin{itemize}
\item If $t^0$ is a nondegenerate critical point of $\Phi_{k,n}(t;z,m)\,$,
then the vector 
$$
v(t^0, z)\,=\, \sum_J A_J(t^0,z) f_Jv
$$
belongs to the subspace ${\rm Sing}(L^{\otimes m})_k=
\{v\in L^{\otimes m} [l(m)-2k]\ \vert\ ev=0\}$
of singular vectors of the weight $l(m)-2k$
and is an eigenvector of operators $H_i(z),\, i=1,...,n$.
\item If $t^0$ is a nondegenerate critical point, then
$$
S(v(t^0,z), v(t^0,z)) = \text{det}_{{}\,1\,\leq \, i,j\leq k\,{}\,}  
\left({\partial^2 \over \partial t_i \partial t_j}\, 
\text{ln}\,\Phi_{k,n}(t^0;z,m)\right)\,.
$$
\item If $l(m) -2k\geq 0$, then for a
generic $z$, the eigenvectors $v(t^0,z)$ generate the space
${\rm Sing}(L^{\otimes m})_k$.
\end{itemize}
\end{thm}

The vectors $v(t^0,z)$ are the Bethe ansatz vectors of the Gaudin model. 
The Bethe equations of the Gaudin model
are the critical point equations for the function $\Phi_{k,n}(t;z,m)$.

\medskip\noindent
{\bf Corollary of Theorems \ref{TM} and \ref{ReV}.}
{\it For a generic 
$z$, the set of the Bethe vectors is a basis in the space ${\rm Sing}(L^{\otimes m})_k$,
that is each eigenvector is presented exactly once as a Bethe vector.}

\medskip\noindent 
Remarks on the Bethe ansatz for $sl(n+1)$ and critical points see in \cite{MV}.

\medskip\noindent
On the connections between the Bethe ansatz and Fuchsian differential equations see
\cite{ Sk1, Sk2}.

\subsection{The function $\Phi_{k,n}(t;z,m)$ and hypergeometric 
solutions of the KZ equations}

The KZ equations \cite {KZ} of the conformal field theory 
for a function $u(z_1,...,z_n)$ with values in the tensor product
$L^{\otimes m}$ is the system of equations
$$
\kappa \,{\partial u\over \partial z_i}\,=\, H_i(z)\,u\,,\qquad i=1,...,n\,.
$$
Here $H_i(z)$ are the operators defined in (\ref{H}). 
The number $\kappa$ is a parameter of equations.

The KZ equations have hypergeometric solutions \cite{SV},
$$
u(z)\,=\, \sum_J \,\int _{\gamma(z)}\, 
\Phi_{k,n}(t;z,m)^{{1\over\kappa}} \,A_J(t,z)
\, dt_1\wedge ... \wedge dt_k\,{}\, f_Jv\,.
$$
The hypergeometric solutions are labeled by suitable families of 
$k$-dimensional cycles $\gamma(z)$. Such a solution takes values 
in the subspace ${\rm Sing}(L^{\otimes m})_k$ of singular vectors.

Studying the semiclassical asymptotics of the hypergeometric 
solutions as $\kappa$ tends to zero
one gets the Bethe ansatz for the Gaudin model \cite{RV, V}.

We see that the function $\Phi_{k,n}(t;z,m)$ is the {\it ``master function''} which governs
solutions of KZ equations and Bethe vectors in the Gaudin model.

\section {Isolated critical points}

\subsection{Combinatorial remarks}
For an $sl_2$ module $V$ and $ l\in \C\,$, let $V[l]=\{v\in V\ \vert\ hv=lv\}$
be the weight subspace of weight $l$. 

Let  $m_1,\dots,m_a\in\ZZ$ and $m_{a+1},\dots,m_n\notin \ZZ$.
Consider the tensor product of irreducible $sl_2$ representations
$
L^{\otimes m}=L_{m_1}\otimes\cdots \otimes L_{m_a}
\otimes L_{m_{a+1}}\otimes\cdots \otimes L_{m_n} \,
$
and for a nonnegative integer $k$ the difference of dimensions
$$
d(m,k)=\dim L^{\otimes m}[l(m)-2k]-\dim L^{\otimes m}[l(m)-2k+2].
$$

Define the number 
$$
\sharp(k,n;m_1,\dots,m_a)=\sum_{q=0}^{a} (-1)^q
\sum_{1\leq i_1<\dots <i_q\leq a}
{k+n-2-m_{i_1}-\dots -m_{i_q}-q\choose n-2}\,.
$$

\medskip\noindent
{\bf Remark.} The formula for $\sharp(k,n;m_1,\dots,m_a)$ 
implies that if $k\leq m_l$ for all $1\leq l\leq a\,$,
then this number does not depend on $m_1,\dots,m_a\,$,
and is equal to ${k+n-2\choose n-2}$.

\begin{thm}\label{RT} We have
$d(m,k)\, =\, \sharp(k,n;m_1,\dots,m_a)\,$.
\end{thm} 

\medskip\noindent
{\bf Proof.}\ \ A basis of the weight subspace 
$L^{\otimes m}[l(m)-2k]$ is formed by the vectors
$
f_Jv=f^{j_1}v_{m_1}\otimes ... \otimes f^{j_n}v_{m_n} \,,
$
where  $j_1, \dots, j_n$ are 
integers such that $j_1+...+j_n=k\,$,  
and for any $l\leq a$ we have
$0\leq j_l\leq m_l$. The Inclusion-Exclusion Principle 
(see \cite{B}, Ch.5) says that the number of such $J$ is equal to
$$
\sum_{q=0}^a\sum_{1\leq i_1<\dots <i_q\leq a} (-1)^q
{k-m_{i_1}-\dots -m_{i_q}-q+n-1\choose n-1} \,.
$$
The Pascal triangle property,
$$
{a\choose b} - {a-1\choose b} = {a-1\choose b-1} \,,
$$
implies the statement. \hfill $\triangleleft$

%{\bf Remark.} 
%If $m_1, ..., m_n \in \Z_{>0},\, k\in \Z_{\geq 0}$,
%then $d(m, l(m)+1-k) \,= \, -d(m,k)\,$.

%\medskip\noindent
%Indeed, if a set $ j_1,...,j_n \in \Z$ satisfies 
%conditions $ 0\leq j_l\leq m_l\,,\  j_1+...+j_n=k$,
%then the set $i_1,...,i_n$ with
%$i_l=m_l-j_l$\, satisfies conditions
%$ 0\leq i_l\leq m_l,\ i_1+...+i_n = l(m)-k$.
%This implies the statement. \hfill $\triangleleft$   

\medskip\noindent
Let $m=(m_1,\dots,m_n)\in \mathbb R^n$, $ k\in\ZZ\,$.

\medskip\noindent
{\bf Definition.} The pair $ \{m, k\}$ is called {\it good} if
$l(m)\geq 2k$ and $m$ has the following form,
$$
m=(m_1,\dots,m_a,m_{a+1},\dots,m_{a+b},m_{a+b+1},\dots,m_n),\ \
0\leq a\leq a+b\leq n \,,
$$
where
\begin{itemize}
\item $m_1,\dots,m_a$ are positive integers;

\item  $m_{a+1},\dots,m_{a+b}$ are positive 
numbers such that for any $1\leq i\leq j\leq b$\\
 the sum $m_{a+i}+\cdots+m_{a+j}$ is not an integer;

\item $m_{a+b+1},\dots,m_n$ are negative integers.
\end{itemize}

\medskip\noindent
{\bf Example.} If $m_1,\dots,m_n, k\in \ZZ\,$, 
and $l(m)\geq 2k\,$, then the pair $\{m,k\}$ is good.

\medskip\noindent
{\bf Remarks.}\ \

\medskip
1. Let the pair  $\{m,k\}$ be good and let 
$
{\rm Sing}(L^{\otimes m})_k=
\{v\in L^{\otimes m} [l(m)-2k]\ \vert\ ev=0\}
$
be the subspace of singular vectors of the weight $l(m)-2k\,$.      
Then $\dim {\rm Sing}(L^{\otimes m})_k\, =\, d(m,k)$.

\medskip
2. If $m_1, ..., m_n,\, k$ are positive integers
and $l(m)-2k \geq 0$, then $w(m,k) = 
\dim {\rm Sing}(L^{\otimes m})_k$.

\medskip\noindent
\begin{lemma} \label{ind}   
Let $p$ be a positive integer, $p<k\,$.
If the pair  $\{(m_1,\dots, m_{j-1},$
$p-1,$ $m_{j+1},\dots,m_n),
k\}$ is good, then the pair 
$\{(m_1,\dots, m_{j-1}, m_{j+1},\dots,m_n, -p-1),
k-p\}$ is good.  \hfill $\triangleleft$  
\end{lemma}

\subsection{Main statements on isolated critical points}
\begin{thm}\label{R}
Let $p$ be a positive integer, $p\leq k\,$.
Assume that the function  $\Phi_{k,n}(t;z,m)$ has an infinite sequence
of critical points  such that
each of its first $p$ coordinates tends to infinity and each of its
remaining coordinates has a finite limit. Then $l(m)=2k-p-1\,$.
\end{thm}
 
\begin{corollary}\label{r} If $l(m) - 2k > -2\,$, then the function
$\Phi_{k,n}(t;z,m)$ has only isolated critical points.\hfill $\triangleleft$
\end{corollary}

\medskip\noindent
{\bf Proof of Theorem \ref{R}.}\ \
We write the system defining the critical
points of  $\Phi(t)=\Phi_{k,n}(t;z,m)$ in the form
$$
\quad (t_r-z_1) (\partial \Phi/\partial t_r)/\Phi=0 \,,
\quad r=1,\dots,k \,.
$$
The $r$-th equation is
$$
-m_1-\sum_{l=2}^n\frac{m_l(t_r-z_1)}{t_r-z_l}+\sum_{1\leq j\leq k
\atop j\neq r}
\frac{2(t_r-z_1)}{t_r-t_j}=0 \,,
$$
and the sum of the first $p$ equations is
$$
-pm_1-\sum_{r=1}^p\sum_{l=2}^n\frac{m_l(t_r-z_1)}{t_r-z_l}
+2\cdot\frac{p(p-1)}2+\sum_{r=1}^p
\sum_{j=p+1}^k\frac{2(t_r-z_1)}{t_r-t_j}=0 \,.
$$
Let $\{ t^{(q)}=(t^{(q)}_1,...,t^{(q)}_k) \}$ be our sequence of
critical points. Then
$$
\frac{t^{(q)}_r-z_1}{t^{(q)}_r-z_l}\underset{q\rightarrow \infty}
{\longrightarrow}\ 1 \,,\ \
\frac{t^{(q)}_r-z_1}{t^{(q)}_r-t^{(q)}_j}\underset{q\rightarrow \infty}
{\longrightarrow}\ 1 \,, \ \ 1\leq r\leq p \,,
\ \ 2\leq l\leq n \,,\ \ p+1\leq j\leq k\,,
$$
and this equation results in $-p(m_1+\dots+m_n)+p(p-1)+2p(k-p)=0\,$.
\hfill $\triangleleft$

\begin{thm}\label{zs}
Let the pair $\{m, k\}$ be good and let $a$ be a nonnegative integer
such that $m_1,\dots,m_a\in\ZZ$ and $m_{a+1},\dots,m_n\notin \ZZ\,$.
%Let  $z^{(s)}=(s,s^2,\dots,s^n)\,$, where $s$ is a real number,
%$s \gg 1\,$.
Then for a generic $z$ in $\C^n$,
all critical points of the function $\Phi_{k,n}(t;z,m)$
are nondegenerate and the critical set consists of
$\sharp(k,n;m_1,\dots,m_a)$ orbits.
\end{thm}     

Theorem \ref{zs} is proved in Sec. ~\ref{Pzs}.  Theorem \ref{zs} implies 
part 1 of Theorem \ref{TM}.

%\begin{thm}\label{NI}
%Let $m_1,\dots,m_n $ and $k$ be fixed positive integers
%such that $l(m)<2k\,$. Let $z^{(s)}=(s,s^2,\dots,s^n)\,$, where
%$s$ is a real number, $s\gg1\,$.
%Then the function $\Phi_{k,n}(t;z^{(s)},m)$ has only non-isolated
%critical points.
%Written in the standard symmetric coordinates $\lambda_1,...,\lambda_k\,$,
%the critical set forms $-d(m,k)$ straight lines
%in the space $\C^k_\lambda\,$.
%\end{thm}     

%Theorem \ref{NI} is proved in Sec. ~\ref{NICP}. 

%Theorem \ref{TM} is a  corollary of Theorems \ref{zs} and \ref{NI}.

%\section {Isolated critical points}\label{ICP}

\subsection{The bound from below}\label{bound}
\begin{thm}\label{ZS}
Let $\{m=(m_1,...,m_n), k\}$ be a good pair. 
Assume that the number $a$ is such that
$m_1,\dots, m_a\in\Z_{>0} \,, \  m_{a+1},\dots,m_n\notin \Z_{>0} \,\,$.
Let $s$ be a real number, $s \gg 1\,$, and  $z^{(s)}=(s,s^2,\dots,s^n)\,$.
Then  the function $\Phi_{k,n}(t;z^{(s)},m)$ has
at least $\sharp(k,n;m_1,\dots, m_a)$ orbits of 
nondegenerate critical points.
\end{thm}         

This Theorem is a direct corollary of results in \cite{RV}, Sec. ~9.
For convenience, we sketch its proof here.

\medskip\noindent {\bf Definition.}\ \
Let $m_1,m_2 \in \R\,$, $k \in \Z_{\geq 0} \,$.
The triple $\{m_1,m_2;k\}$ is called {\it admissible} 
if the following two conditions are satisfied,
\begin{itemize}
\item  $m_1+m_2 - 2k \geq 0$\,,

\item  if for some $i \in \{1,2\}$ we have $m_i\in \Z_{\geq 0}\,$, then 
$k \leq m_i$\,.
\end{itemize}
If the triple $\{m_1,m_2;k\}$ is admissible, then the function
$$
\Phi_{k, 2}(t) = \prod_{i=1}^{k}
t_i^{-m_1}(t_i-1)^{-m_2}\prod_{1\leq i<j\leq k}(t_i-t_j)^2
$$
has exactly $k!$  critical points all of which are nondegenerate
\cite{V}.

\medskip\noindent {\bf Definition.}\ \
Let $\{m=(m_1,...,m_n),k\}$ be a good pair. 
Let $I=(i_1,\dots,i_n)$ be a sequence of nonnegative integers 
such that $i_1=0$ and $i_2+\cdots+i_n=k\,$. The sequence $I$ is called
{\it an admissible sequence} for $\{m,k\}$ if all triples
$$
\{m_1+\cdots+m_{l-1}-2(i_1+\cdots+i_{l-1}), \,m_{l}; \,i_{l}\}\,,
$$
for $  l=2,\dots,k$
are admissible.

\medskip\noindent 
{\bf Proof of Theorem \ref{ZS}.}\ \ 
For an admissible sequence $I=(i_1,\dots,i_n)\,$, make a change of
variables
$$
t_j=s^lu_j\ \ {\rm if}\ \ i_1+\cdots+i_{l-1}<j\leq i_1+\cdots+i_l \,,
\ \ l=2,\dots,n \,.
$$
For $2\leq l\leq n\,$, define the function
\begin{eqnarray*}
\Phi_{i_l,2} & = & \Phi_{i_l,2}(u_{i_1+\dots+i_{l-1}+1},\dots,
u_{i_1+\dots+i_l})\\
 & = & \prod_{j=i_1+\dots+i_{l-1}+1}^{i_1+\dots+i_l}
u_j^{-a_l}(u_j-1)^{-m_l}\prod_{i_1+\dots+i_{l-1}+1\leq i<j\leq
i_1+\dots+i_l}(u_i-u_j)^2 \,,
\end{eqnarray*}  
where $a_l =m_1+\cdots+m_{l-1}-2(i_1+\cdots+i_{l-1})\,$.

Let $\Phi_I(u)=\Phi_{i_2,2}\cdots\Phi_{i_n,2}\,$.
For any $l=2,\dots,n\,$, the function $\Phi_{i_l,2}$
has exactly one orbit of nondegenerate critical points according to 
Theorem 1.3.1 in \cite{V}.
Let $u_{(l)}$ be a critical point of $\Phi_{i_l,2}\,$, then
$u_I=(u_{(2)},\dots,u_{(n)})$ is a nondegenerate critical point
of the function $\Phi_I(u)\,$. In a neighborhood of $u_I\,$,
the critical point system of the function 
$\Phi_{k,n}(t(u))=\Phi_{k,n}(t(u);z^{(s)},m)$ 
is a deformation of the critical point system of the function
$\Phi_I(u)$ with deformation parameter $s\,$,
$$
\frac{\partial \Phi_{k,n}(t(u))/\partial u_j}{\Phi_{k,n}(t(u))}=
\frac{\partial \Phi_I(u)/\partial u_j}{ \Phi_I(u)} + O(s^{-1})=0 \,,
\ \ j=1,\dots,k \,.
$$
When $s \to \infty\,$, the function $\Phi_{k,n}(t(u))$ 
has a nondegenerate critical point $u_I(s)$ close to $u_I\,$,
which defines a nondegenerate critical point $t_I(s)$ of 
the function $\Phi_{k,n}(t;z^{(s)},m)\,$. Theorem 9.9 in \cite{RV}
and its corollaries imply that if $I$ and $I'$ are distinct
admissible sequences,
then the corresponding points $t_I(s)$ and $t_{I'}(s)$ cannot
belong to the same orbit. To complete the proof it remains to note
that the number of admissible sequences for $\{m, k\}$ is equal to
the dimension of  ${\rm Sing}(L^{\otimes m})_k \,$. 
This follows from the fact that the admissible sequences label 
a basis of iterated singular vectors in 
${\rm Sing}(L^{\otimes m})_k \,$, see Sec. 8 in \cite{RV}.
\hfill $\triangleleft$

\subsection{The maximal possible number of critical points}

\begin{thm}\label{k1} 
%Let $z^{(s)}=(s,s^2,\dots,s^n)\,$, where 
%$s\in\R\,$, $s\gg 1\,$.
If the pair $\{m,k\}$ is good and if $k\leq m_i$ for
all $m_i\in\Z_{>0}\,$, then for a generic $z$ in $\C^n$,
the function $\Phi_{k,n}(t; z, m)$
has exactly ${k+n-2\choose n-2}$ orbits of critical points which all
are nondegenerate.
\end{thm}
 
\begin{corollary}\label{k=1}
Theorem \ref{zs} is true for $k=1\,$.
\end{corollary}            
\medskip\noindent 
{\bf Proof of Theorem \ref{k1}.}\ \ If all 
numbers $z_1,...,z_n$ are real and all numbers $m_1,..., m_n$
are negative, then all critical points of the function
$\Phi_{k,n}(t;z,m )\,$
are nondegenerate and the critical set consists of
${k+n-2\choose n-2}$ orbits \cite{V}. Therefore for any $z$ and $m$
the total number of isolated orbits of critical points counted
with multiplicities is not greater than ${k+n-2\choose n-2}\,$.
Corollary \ref{r} says that the function $\Phi_{k,n}(t; z, m )$
does not have non-isolated critical points. In order to finish the proof
we apply Theorem \ref{ZS}. \hfill $\triangleleft$                   

\subsection{Vanishing critical points}  \ \ 

\medskip\noindent 
Set $m(\epsilon)=(m_1,\dots, m_{j-1}, k-1+\epsilon,m_{j+1},\dots,m_n)\,$.
Let $K$ be the number of critical points of the function 
$\Phi_{k,n}(t;z,m(\epsilon))$ which tend to the vertex
$\{t_1=\dots=t_k=z_j\}$ when $\epsilon$ tends to zero. 

\begin{thm}\label{k}
If the pair   $\{(m_1,\dots, m_{j-1},k-1,m_{j+1},\dots,m_n),k\}$
is good, then the number $K$ is positive and divisible by $k!$\, \ .
\end{thm}
             
Theorem \ref{k} is proved in Sec. \ref{Pk}.

\medskip
Let $p$ be a positive integer, $p<k\,$. For
 $m=(m_1,\dots, m_{j-1},p-1,m_{j+1},\dots,m_n)\,$,
 set $m(\epsilon)=(m_1,\dots, m_{j-1},
p-1+\epsilon,m_{j+1},\dots,m_n)\,$,
$m^{(p)}=(m_1,...,m_{j-1},-p-1,m_{j+1},\dots,m_n)\,$.

\medskip\noindent
{\bf Definition.}\ \ The function
$$
\Phi_{k-p,n}(t_{p+1},...,t_k;z,m^{(p)})
 = \prod_{i=p+1}^k(t_i-z_j)^{p+1}
\prod_{i=p+1}^k\prod_{l\neq j}(t_i-z_l)^{-m_l}
\prod_{p+1\leq i<j\leq k}(t_i-t_j)^2 \,
$$
is  called the function  {\it induced} by the function
$\Phi_{k,n}(t;z,m(\epsilon))$
on the edge $\{t_1=\dots=t_p=z_j\}$ as $\epsilon$ tends to zero.

\medskip
Let $\epsilon$ tend to zero.
Let  $B=(b_{p+1},\dots,b_k)$ be a nondegenerate critical point of
the induced function $\Phi_{k-p,n}(t_{p+1},...,t_k;z,m^{(p)})\,$, 
and let $K$ be the number of critical points of the function
$\Phi_{k,n}(t;z,m(\epsilon))$ which tend to the point
$\{t_1=\dots=t_p=z_j, t_{p+1}=b_{p+1},\dots,t_k=b_k\}\,$.     

\begin{thm}\label{p} 
If $\{(m_1,\dots, m_{j-1},p-1,m_{j+1},\dots,m_n),k\}$
is a good pair, then $K$ is positive and divisible by $p!$\, \ . 
\end{thm}

Theorem \ref{p} is proved in Sec. \ref{Pp}.

\subsection{Proof of Theorem \ref{k}}\label{Pk}

\medskip\noindent 
After the translation $t_i\mapsto t_i-z_j\,$, $z_l\mapsto z_l-z_j\,,$  
and renumbering $z_1,\dots, z_n\,$, we can assume
$z=(0,z_2,\dots,z_n)$ and
$m(\epsilon)=(k-1+\epsilon,m_2,\dots, m_n)\,$.
We estimate the number of critical points $t(\epsilon)$
of the function
$$
\Phi_{k,n}(t;z,m(\epsilon))=\prod_{i=1}^k\left[t_i^{-k+1-\epsilon}
\prod_{l=2}^n(t_i-z_l)^{-m_l}\right]
\prod_{1\leq i<j\leq k}(t_i-t_j)^2 \,
$$
such that $t_r(\epsilon)$ tends to zero for $ r=1,\dots,k$ 
as $\epsilon$ tends to zero.
 
Blow-up  the vertex $\{t_1=\dots=t_k=0\}\,$.
In  coordinates
$u_1\,, \dots\,, u_k\,$, where
$$
t_1=u_1u_k,\ \ \dots,\ \ t_{k-1}=u_{k-1}u_k,\ \
t_k=u_k,
$$
the function $\Phi_{k,n}(t;z,m(\epsilon))$ has the  form

\begin{eqnarray*}
\tilde\Phi & = &
\Phi_{k,n}(u_1u_k,\dots, u_{k-1}u_k, u_k;z,m(\epsilon))= \\
 & = & \prod_{i=1}^{k-1} u_i^{-k+1-\epsilon}
(u_i-1)^2\prod_{1\leq i<j\leq k-1}(u_i-u_j)^2 \\
 & \cdot  &  u_k^{-k\epsilon}\prod_{l=2}^n(u_k-z_l)^{-m_l}
\prod_{i=1}^{k-1}\prod_{l=2}^{n}(u_iu_k-z_l)^{-m_l} \,.
\end{eqnarray*}     

Consider this function as a function on the space $\C^{k+1}$ with
coordinates $u_1\,,\dots\,,u_k\,,\ \epsilon\,$. Consider in $\C^{k+1}$
the set $C$ of all critical points of this function with respect to
coordinates $u_1,\dots,u_k\,$.
 
\begin{lemma}\label{u_k} 
Near the divisor $\{ \, (u; \epsilon)\in \C^{k+1} \, | \, u_k=0\,\}\,$, 
the critical set $C$ is the union of $(k-1)!$ nonsingular curves
which intersect the divisor at $(k-1)!$ points
$(A_{\sigma},0;0)\,$, where $A_{\sigma}$ runs through all
permutations of $\{\alpha,\dots,\alpha^{k-1}\},\
\alpha=\exp(2\pi i/k)\,$.
The coordinate $u_k$ is a local parameter at the intersection point on
each of these curves.
\end{lemma}
\medskip\noindent 
{\bf Proof.}\ \
We write the system defining the critical points of $\tilde\Phi$
in the form
$$
\frac{\partial \tilde\Phi/\partial u_q}{\tilde\Phi}=0 \,,
 \ q=1,\dots, k-1,\ \quad
u_k\cdot\frac{\partial \tilde\Phi/\partial u_k}{\tilde\Phi}=0 \,,
$$
and get the following equations   

\begin{eqnarray*}
&& \frac{-k+1-\epsilon}{u_q}+\frac{2}{u_q-1}+\sum_{{1\leq j\leq k-1}
\atop{j\neq q}}\frac{2}{u_q-u_j}
-\sum_{l=2}^nm_l\frac{u_k}{u_qu_k-z_l}=0 \,, \\
&& -k\epsilon-\sum_{l=2}^nm_l\frac{u_k}{u_k-z_l}
-\sum_{i=1}^{k-1}\sum_{l=2}^{n}m_l\frac{u_iu_k}{u_iu_k-z_l}=0 \,,
\end{eqnarray*}
where $ q=1,\dots,k-1\,$.
 
One can express $\epsilon$ in terms of $u_1,\dots,u_k$ from
the last equation, and consider the  equations
$$
\frac{\partial \tilde \Phi/\partial u_q}{\tilde \Phi} =0 \,,
\ \ q=1,\dots,k-1,
$$
as a system of equations  with respect to $u_1,\dots, u_{k-1}$
depending on the parameter $u_k\,$.
For $u_k=0\,$, this system turns into the critical point system
of the function
$$
\Phi_{k-1,2}(u)=\prod_{i=1}^{k-1}u_i^{-k+1}(u_i-1)^2
\prod_{1\leq i<j\leq k-1}(u_i-u_j)^2 \,.
$$
Theorem 1.3.1 \cite{V} implies that the function
$\Phi_{k-1,2}$ has exactly $(k-1)!$ critical points,
all of which are non-degenerate, and the coordinates of
each of these critical points form the set of all roots of the 
equation $\xi^{k-1}+\xi^{k-2}+\dots+1=0\,$.
This gives the Lemma. \hfill $\triangleleft$
 
\begin{lemma}\label{AV} For a given permutation $A_\sigma\,$,
the number of critical points     
of the function
$\tilde\Phi\,|_{\epsilon =\epsilon_0}\,$, which tend to
$(A_{\sigma},0;0)$ as
$\epsilon_0$ tends to zero,
is positive and divisible by $k\,$.
\end{lemma}
\medskip\noindent 
{\bf Proof.}\ \ It is enough to prove the statement for
$A=(\alpha,\alpha^2,\cdots,\alpha^{k-1})\,$.
The function $\Phi_{k,n}(t)$ is invariant with respect to
permutations of $\{t_1,\dots,t_k\}\,$. Therefore
the critical set $C$ is invariant with respect
to the corresponding action of the symmetric group $S^k$
on the space $\C^{k+1}$ with coordinates 
$u_1\,,\cdots\,, u_k\,,\ \epsilon\,$.
The connected component $C_A \subset C$ which contains  $(A,0;0)$
is preserved by the map $\PP\,$, the lifting of the cyclic permutation
$ t_1 \mapsto t_2 \mapsto \dots \mapsto t_k \mapsto t_1$; and
the point $(A,0;0)$ is a fixed point of $\PP\,$. According to
Lemma \ref{u_k}, the coordinate $u_k$ is a local parameter
on the curve $C_A\,$,
$$
C_A=\{ \, (u;\epsilon)\ \vert \ u_j=\alpha^j+ O(u_k),\
j=1,\dots, k-1,\ \  \epsilon=f(u_k)\, \} \,,
$$
where $f(u_k)$ is the germ of a suitable holomorphic function.
This germ can not be identically zero, as in this case the 
function $\Phi_{k,n}(t;z,m)$ would have a curve of critical
points, 
$$
\{\, t_j=\alpha^jt_k+O(t_k)\,,\quad  j=1,\dots,k-1\, \}\,,
$$
but this is impossible by Corollary \ref{r}
because the pair $\{(m_1,\dots, m_{j-1},k-1,m_{j+1},\dots,m_n),k\}$ 
is good. The equation $\epsilon = f(u_k)$ has to be
invariant with respect to the map $\PP$ which does not change
$\epsilon$ and maps $u_k$ to
$ u_1u_k=\alpha u_k + O(u_k^2) \,$.
This means that the Taylor expansion of the germ $f(u_k)$
starts with a power of $u_k$ divisible by $k\,$.   \hfill $\triangleleft$
 
Lemmas \ref{u_k} and \ref{AV} imply Theorem \ref{k}.  \hfill $\triangleleft$    

\subsection{Proof of Theorem \ref{p}}\label{Pp}
We set 
$$
z=(0,z_2,\dots, z_n)\,,\quad 
m(\epsilon)=(p-1+\epsilon, m_2,\dots, m_n)\,,
$$
and count the number of critical points $t(\epsilon)$ of the function
$$
\Phi_{k,n}(t;z,m(\epsilon))=
\prod_{i=1}^k\left[t_i^{-p+1-\epsilon}\prod_{l=2}^n
(t_i-z_l)^{-m_l}\right]\prod_{1\leq i<j\leq k}(t_i-t_j)^2
$$
which satisfy
$$
t_i(\epsilon)\rightarrow\, 0 \,,
\ \ 1\leq i\leq p, \quad
t_j(\epsilon)\rightarrow\, b_j \,,\ \   p+1\leq j\leq k,
$$
as $\epsilon$ tends to $0\,$.

Blow-up the edge $\{t_1=\dots=t_p=0\}\,$. In coordinates
$u=(u_1,\dots,u_p)\,$, $t'=(t_{p+1},\dots,t_k)\,$, where
$t_1=u_1u_p,\ \dots,\ t_{p-1}=u_{p-1}u_p,\ t_p=u_p\,,$
the function $\Phi_{k,n}(t;z,m(\epsilon))$ is
\begin{eqnarray*}
 \tilde\Phi & = & \Phi_{k,n}(u_1u_p,\dots,u_{p-1}u_p, u_p,
t_{p+1},\dots,t_k;z,m(\epsilon))=\\
& & \prod_{i=p+1}^k\left[t_i^{-p+1-\epsilon}
\prod_{l=2}^n(t_i-z_l)^{-m_l}\right]
\prod_{p+1\leq i<j\leq k}(t_i-t_j)^2 \\
& & \cdot  \prod_{i=1}^{p-1}\prod_{j=p+1}^k(t_j-u_iu_p)^2
    \prod_{j=p+1}^k(t_j-u_p)^2 \\
& & \cdot \prod_{i=1}^{p-1} u_i^{-p+1-\epsilon} 
(u_i-1)^2\prod_{1\leq i<j\leq p-1}(u_i-u_j)^2 \\
& & \cdot u_p^{-p\epsilon}\prod_{l=2}^n(u_p-z_l)^{-m_l}
\prod_{i=1}^{p-1}\prod_{l=2}^n(u_iu_p-z_l)^{-m_l} \,.
\end{eqnarray*}

We take the critical point system for $\tilde\Phi$ 
in the following form
\begin{eqnarray*}
 \frac{\partial\tilde\Phi/\partial u_i}{\tilde\Phi}=0 \,,
\ \ & i=1,\dots,p-1; & \quad (S_u)\\
 \frac{\partial\tilde\Phi/\partial t_j}{\tilde\Phi}=0 \,,
\ \ & j=p+1,\dots,k; &  \quad (S_{t'})\\
 u_p\cdot\frac{\partial\tilde\Phi/\partial u_p}{\tilde\Phi}=0 \,. 
& & \quad (S_p)
\end{eqnarray*}

From  equation $(S_p)\,$, one can express $\epsilon$ in terms of
$u\,,t'\,$. Therefore one can consider equations
$(S_u),\ (S_{t'})$ as a system
of equations  with respect to $u_1,\dots, u_{p-1}\,$,
$t_{p+1},\dots, t_k$ depending on the parameter $u_p\,$. For $u_p=0\,$,
equations $(S_u)$ turn into the critical point system of the
function
$$
\Phi_{p-1,2}(u)=\prod_{i=1}^{p-1}u_i^{-p+1}(u_i-1)^2
\prod_{1\leq i<j\leq p-1}(u_i-u_j)^2 \,, 
$$
and equations $(S_{t'})$ turn into the
critical point system of the induced function
$$
\Phi_{k-p,n}(t';z,m^{(p)}) =
\prod_{i=p+1}^k\left[t_i^{p+1}\prod_{l=2}^n(t_i-z_l)^{-m_l}\right]
\prod_{p+1\leq i<j\leq k}(t_i-t_j)^2 \,.
$$

Consider $\tilde\Phi$ as a function on the space $\C^{k+1}$
with coordinates $u\,,t'\,,\epsilon\,$. Consider in $\C^{k+1}$ the
critical set of $\tilde\Phi$ with respect to $u\,,t'\,$.
Similarly to Lemma \ref{u_k}  we get

\begin{lemma}\label{u_p}
The critical  set near the plane
$\{\, (u,t',\epsilon)\in\C^{k+1}\, \vert\, 
u_p=0,\, t'=B\, \}$ is the union of $(p-1)!$
nonsingular curves which intersect this  plane 
at $(p-1)!$ points $(A_{\sigma},0,B;0)\,$, where $A_{\sigma}$ 
runs through all permutations of $\{\alpha,...,\alpha^{p-1}\},\
\alpha=\exp(2\pi i/p)\,$.
The coordinate $u_p$ is a local parameter 
at the intersection point on each these curves.\hfill $\triangleleft$
\end{lemma} 
The function $\Phi_{k,n}(t)$ is invariant with respect to
permutations of $\{ t_1,\dots,t_p \}\,$, and hence
the union of these $(p-1)!$ curves is invariant with respect 
to the corresponding action of the symmetric group $S^p$ on
the space $\C^{k+1}$ with coordinates $u\,,t',\ \epsilon\,$.

\begin{lemma}\label{av}
For a given permutation $A_{\sigma}\,$, the number of critical
points of the function $\tilde\Phi |_{\epsilon=\epsilon_0}$
which tend to $(A_{\sigma}, 0, B; 0)$ as $\epsilon_0 $ tends 
to zero is positive and divisible by $p\,$.
\end{lemma}

\medskip\noindent  
{\bf Proof.}\ \ We prove this statement for $A=(\alpha,\dots,
\alpha^{p-1})\,$.
The connected component of the critical set which contains
the point $(A,0,B;0)$ is of the form
$$
C_{A,B}=\{\, u_i=\alpha^i+O(u_p),\, i=1,\dots,p-1;\ \ t_j=b_j+O(u_p),\,
j=p+1,\dots,k;\ \ \epsilon=f(u_p)\, \} \,,
$$
where $f(u_p)$ is the germ of a suitable holomorphic function.
Similarly to Lemma \ref{AV}, we conclude that $f$ is a non-zero
germ, that $C_{A,B}$ is invariant with respect to the map $\PP$ which 
is the lifting of the permutation
$ t_1\mapsto t_2\mapsto\dots\mapsto t_p\mapsto t_1\,$, and
that the Taylor
expansion of the germ $f(u_p)$ starts with a power of $u_p$ 
divisible by $p\,$. \hfill $\triangleleft$

Lemmas \ref{u_p} and \ref{av} imply Theorem \ref{p}.  \hfill $\triangleleft$    

\subsection{Proof of Theorem \ref{zs}}\label{Pzs}

\medskip\noindent 
We prove the statement by a double induction with respect to $k\,$, 
the number of variables in $\Phi_{k,n}(t; z, m)\,$, and $a(m)\,$, 
the number of positive integers in $m\,$.
 
For $k=1$ and any $a(m)$ the statement is true by Corollary
\ref{k=1}. For any $k$ and $a(m)=0\,$, the statement holds by Theorem
\ref{k1}.
 
Assume that the Theorem is proved for $k < k_0$ and any $a(m)$ and 
for $k = k_0$ and $a(m) < a_0\,$. We prove the Theorem for $k=k_0$ 
and $a=a_0\,$.
 
Let $\{m=(m_1,\dots,m_n),k\}$ be a good pair.
Assume that $m_1,\dots,m_a\in\ZZ$ and
$m_{a+1},\dots,m_n\notin\ZZ$.
For $\epsilon\neq 0$ small enough, the pair
$\{m(\epsilon)=(m_1,\dots,m_{a-1},m_a+\epsilon,m_{a+1},\dots,m_n),k\}$
is also good, and the number of positive integers in $m(\epsilon)$
is $a-1\,$. Therefore according to the induction hypothesis, for a generic $z$ the
function $\Phi_{k,n}(t;z, m(\epsilon))$  has exactly
$\sharp(k,n;m_1,\dots,m_{a-1})$ orbits of critical points which all
are nondegenerate.
 
We study how the number of orbits of critical points of 
$\Phi_{k,n}(t;z, m(\epsilon))$ changes as 
$\epsilon\rightarrow 0\,$.
According to Corollary \ref{r}, non-isolated critical points do not
appear. For isolated critical points, there are three possibilities.
 
\medskip\noindent
(1)\ \  If $m_a\geq k\,$, then the function  $\Phi_{k,n}(t;z, m)$
has at most
$$
\sharp(k,n;m_1,\dots,m_{a-1}) = \sharp(k,n;m_1,\dots,m_{a})
$$
orbits of critical points. Indeed, the number of orbits of
isolated critical points cannot increase as $\epsilon\rightarrow 0\,$.
 
\medskip\noindent
(2)\ \  If $m_a=k-1\,$, then according to Theorem \ref{k} at least $k!$
critical points disappear at the vertex $t_1=\dots=t_k=z_a$ as
$\epsilon\rightarrow 0\,$.
Therefore the number of orbits of critical points of the function
$\Phi_{k,n}(t; z, m)$ does not exceed
$$
\sharp(k,n;m_1,\dots,m_{a-1}) - 1 = \sharp(k,n;m_1,\dots,m_{a-1},k-1)
= \sharp(k,n;m_1,\dots,m_a) \,.
$$
     
\medskip\noindent
(3)\ \  If $m_a=p-1$ for some integer $1< p\leq k-1\,$, then critical
points disappear at certain points of the edges of the form
$\{t_{i_1}=\dots=t_{i_{p}}=z_a\}\,$.
These points are critical points of the functions induced by the
function  $\Phi_{k,n}(t;z, m(\epsilon))$ on the edges as
$\epsilon\rightarrow 0\,$. The number of the edges is ${k\choose
k-p}\,$. Any of the induced functions is a function of $k-p<k$
variables with the set of exponents 
$$
m^{(p)}=(m_1,\dots,m_{a-1}, -p-1, m_{a+1},\dots,m_n) \,.
$$
The pair $\{m^{(p)},k-p\}$ is  good after renumbering the coordinates of
the vector $m^{(p)}$,  and $m^{(p)}$ contains $a-1$ positive integers. 
Hence according to the induction hypothesis, for a generic $z$
any of the induced functions has exactly
$(k-p)!\sharp(k-p,n;m_1,\dots,m_{a-1})$ critical points which all
are nondegenerate. At each of these points, at least $p!$
critical points of the function $\Phi_{k,n}(t; z, m(\epsilon))$
disappear as $\epsilon\rightarrow 0$ by Theorem \ref{p}. Thus
the total number of critical points which disappear as 
$\epsilon\rightarrow 0$ is at least
$$
{k\choose k-p} \, (k-p)! \,  \ p! \, \ \sharp(k-p,n;m_1,\dots,m_{a-1}) \,= \,
k!\,\ \sharp(k-m_a-1,n;m_1,\dots,m_{a-1}) \,.
$$
Therefore  $\Phi_{k,n}(t;z^{(s)},m)$ has at most
$$
\sharp(k,n;m_1,\dots,m_{a-1})-\sharp(k-m_a-1,n;m_1,\dots,m_{a-1})=
\sharp(k,n;m_1,\dots,m_{a})
$$
orbits of critical points. 

\medskip
Thus in all cases the number of orbits of critical points of the
function $\Phi_{k,n}(t; z, m)$
is not greater than $\sharp(k,n;m_1,\dots,m_{a})\,$. 
But Theorem  \ref{ZS} says that the number of orbits is at least
$\sharp(k,n;m_1,\dots,m_{a})\,$. This gives Theorem \ref{zs}. \hfill $\triangleleft$                           

\section{Critical points and Fuchsian equations}\label{NICP}

\subsection{Critical points and Fuchsian equations  
with only polynomial solutions}

\medskip\noindent 
On Fuchsian equations see  \cite{R}.

\medskip\noindent 
\begin{lemma}\label{op1} ${}$
Let all solutions of the Fuchsian equation (\ref{Ee}) be polynomials.
Then generic solutions have no multiple roots.
\end{lemma}

\medskip\noindent 
{\bf Proof.}\ \
Let $v(x)$ be a solution. Assume
that the order of $v(x)$ at some point $x=z_0$  is $r$, \, $r\geq 2$. Then
the order of $F(x)v''(x)$ at $x=z_0$ is at least $r-1$. Hence $F(z_0)=0$.
Therefore $z_0$ is one of the points $z_1,...,z_n$ and  the order
of $v(x)$ at this point is $m_j+1$. This means that $v(x)$ is not a
generic solution.  \hfill $\triangleleft$

\medskip\noindent 
\begin{lemma}\label{op2} ${}$
Let $m=(m_1,...,m_n) \in \Z^n_{>0}, \, k\in \Z_{>0}$. Let
$t^0$ be a critical point of the function $\Phi_{k,n}(t;z,m)$.
Then all solutions of the associated differential equation
$E(t^0,z,m)$ are polynomials.
\end{lemma}

\medskip\noindent 
{\bf Proof.}\ \ Let $u(x)=(x-t_1^0)\cdots (x-t^0_k)$. 
For $j=1,...,n$, we have $u(z_j)\neq 0$. Hence
all solutions are univalued at $z_j$. 
Therefore all solutions are univalued at infinity as well. 
Thus all solutions are  polynomials. \hfill $\triangleleft$

\medskip\noindent
{\bf Remarks.}\ \ 

\medskip\noindent 
1. If $l(m)+1-k > k$, then the generic solution of equation
$E(t^0,z,m)$  has degree $l(m)+1-k$.

\medskip\noindent 
2. If $\ 0\leq l(m)+1-k < k $, then 
the generic solution of equation
$E(t^0,z,m)$ has degree $k$, the equation also has
solutions of degree $l(m)+1-k$.

\medskip\noindent 
3. If $l(m)+1-k = k $, then the two exponents at infinity are equal.
Every Fuchsian differential equation with equal exponents has 
multivalued solutions.
Hence the function $\Phi_{k,n}(t;z,m)$ does not have critical points.
This is the second part of Theorem \ref{TM}.

\medskip\noindent 
4. If $l(m)+1-k <0 $, then one of exponents at infinity is positive.
Such a Fuchsian differential equation cannot have only polynomial
solutions. Hence the function $\Phi_{k,n}(t;z,m)$ does not have critical points.
This is the fourth part of Theorem \ref{TM}.

\medskip\noindent 
5. Let $l(m)+1-k=0$ and let equation (\ref{Ee}) have only polynomial solutions
with the degree of the generic solution equal to $k$. Then $H(x)$ is identically equal to zero
and the solutions have the form 
$$
\int (x-z_1)^{m_1}\cdots (x-z_n)^{m_n} dx + \text{const}\,.
$$
Hence the critical set of the function 
$\Phi_{k,n}(t; z, m)$, written in symmetric coordinates,
forms a straight line. In this case  $w(m,l(m)+1-k)=1$.    
This statement gives part 3 of Theorem \ref{TM} for $l(m)+1-k=0$.

\medskip\noindent 
\begin{lemma}\label{op3} ${}$
Let $m=(m_1,...,m_n) \in \Z^n_{>0}, \, k\in \Z_{>0}, \, l(m) -2k\leq -2$. Let
$t^0$ be a critical point of the function $\Phi_{k,n}(t;z,m)$.
Then there exists a curve of critical points containing $t^0$. 
The curve being written in symmetric
coordinates $\lambda_1=\sum t_i, \, ... \,,\, \lambda_k = t_1 \cdots t_k$ 
is a straight line in $\C^k_\lambda$.
\end{lemma}

\medskip\noindent 
{\bf Proof.}\ \ Equation $E(t^0,z,m)$ has only polynomial solutions. Let
$u_1(x)=(x-t_1^0)\cdots (x-t_k^0)$ and let $u_2(x)$ be a solution of degree
$l(m)+1-k$. Then solutions $u_c(x)=u_1(x) + cu_2(x)$ correspond 
to a curve of critical points. The coefficients of $u_c(x)$ give 
a straight line in $\C^k_\lambda$. \hfill $\triangleleft$

\medskip\noindent 
\begin{lemma}\label{op4} ${}$
Let $m=(m_1,...,m_n) \in \Z^n_{>0}, \, k\in \Z_{>0}, \, l(m) -2k\leq -2$.
Then the straight lines in $\C^k_\lambda$ of critical points of 
$\Phi_{k,n}(t;z,m)$ do not intersect.
\end{lemma}

\medskip\noindent 
{\bf Proof.}\ \ If two critical points of 
$\Phi_{k,n}(t;z,m)$ belong to different lines,  then the associated 
differential equations are different. Two differential equations
of the form (\ref{Ee}) with the same $F(x), G(x)$ and distinct $H(x)$
cannot have common nonzero solutions. \hfill $\triangleleft$

\medskip\noindent 

If $k$ is such that $l(m)-2k\leq -2$, then for
 $k'=l(m)+1-k$ we have $l(m)-2k' \geq 0$.

\medskip\noindent 
\begin{lemma}\label{op5} ${}$
Let $m=(m_1,...,m_n) \in \Z^n_{>0}, \, k\in \Z_{>0}$, $l(m) -2k\leq  -2$. 
Then the number of critical lines in $\C^k_\lambda$ 
of the function $\Phi_{k,n}(t;z,m)$ is not less than
the number of orbits of critical points of the function
$\Phi_{l(m)+1-k,\,n}(t;z,m)$.
\end{lemma}

\medskip\noindent 
{\bf Proof.}\ \ Let $t^0\in \C^{l(m)+1-k}$ be a critical point of 
$\Phi_{l(m)+1-k,n}(t;z,m)$. Generic solutions of $E(t^0,z,m)$ are of degree $k$. 
They define a straight line
in $\C^k_\lambda$ of critical points of $\Phi_{k,n}(t;z,m)$.

If two critical points of $\Phi_{l(m)+1-k,\,n}(t;z,m)$ belong to different orbits, 
then the associated differential equations are different. 
The corresponding straight lines do not intersect. \hfill $\triangleleft$

\medskip\noindent 
\begin{thm}\label{Prop} Let $m=(m_1,...,m_n) \in \Z^n_{>0}, \, k\in \Z_{>0}$, 
$0 < l(m)+1-k < k$. 
%Let  $z^{(s)}=(s,\, s^2,\, \dots, \, s^n)$, 
%s\in\R,\,$  $s\gg 1\,$. 
For a generic $z$ in $\C^n$, let $t^0$ be a critical point of the function
$\Phi_{k,n}(t; z, m)$. Let $u(x)$ be a solution of degree $l(m)+1-k$ of equation
$E(t^0, z, m)$. Then roots of $u(x)$ form a critical point of the function
$\Phi_{l(m)+1-k,\,n}(t; z, m)$.
\end{thm}
 
\medskip\noindent 
\begin{corollary}\label{C}
Under conditions of Theorem \ref{Prop}, for a generic $z$
the number of critical lines of the function
$\Phi_{k,n}(t; z, m)$ is equal to the number of 
orbits of critical points of the function $\Phi_{l(m)+1-k,\,n}(t; z, m)$.
\end{corollary}            

\medskip\noindent 
Corollary \ref{C} implies Theorem \ref{main2} and part 3 of Theorem \ref{TM}.

Theorem \ref{Prop} is proved in Sections \ref{mult} and \ref{PNI}. 
To prove Theorem \ref{Prop} one needs  to show that if $z$ is generic, then
the polynomial $u(x)$ does not have multiple roots.

\subsection{Polynomial solutions with multiple roots}\label{mult}
\medskip\noindent
Let $z_1,\dots,z_n$ be pairwise distinct complex numbers. Let
$m_1,\dots,m_n, \, k\in \Z_{>0}$.
Let $F(x)$  and  $G(x)$ be determined by (\ref{M}).
Let $H(x)$ be a polynomial of degree not greater than $n-2$ such that the differential
equation (\ref{Ee}) has a polynomial solution $U(x)$ of degree $k$ with a multiple root.
Then the root is equal to one of the singular points $z_j$ and the multiplicity of the root is
equal to $m_j+1$. 

Assume that $z_{a}, z_{a+1},...,z_n$ are all multiple roots of the solution, where $a$ is 
a  suitable number.  Then
\begin{eqnarray}\label{mr}
U(x) \,=\,  (x-t_1^0)\cdots(x-t_b^0) \, \cdot \, 
(x-z_{a})^{m_{a}+1} \cdots (x-z_n)^{m_n+1}\, 
\end{eqnarray}
where $t_1^0,...,t_b^0$ are roots of multiplicity 1 and $b + (m_{a}+1)
+...+(m_n+1)=k$.

\medskip\noindent 
\begin{lemma}\label{ps}\ \

\begin{itemize}
\item Under the above assumptions,  $t^0=(t_1^0,...,t_b^0)$ is a critical point
of the function $\Phi_{b, n}(t; z, \tilde m)$ where
$\tilde m = (m_1,...,m_{a-1}, -m_{a}-2,..., -m_n-2)$.

\item If  $t^0=(t_1^0,...,t_b^0)$ is a critical point
of the function $\Phi_{b, n}(t; z, \tilde m)$,
then there exists a unique polynomial $H(x)$ of degree not greater than $n-2$
such that the polynomial
$U(x)$ given by  (\ref{mr}) is a polynomial solution of equation (\ref{Ee}).

\end{itemize}
\end{lemma}

%Lemma \ref{ps} is a simple analog of Theorem \ref{PS}.

The differential equation of part 2 of Lemma \ref{ps} will be called {\it associated with
the critical point $t^0$ and vectors $m$, $\tilde m$} and denoted $E ( t^0, z, m, \tilde m ) $.

\medskip\noindent
{\bf Proof.} 
The substitution $x=t_i^0$ into (\ref{Ee}) gives
$$
\frac {U''(t_i^0)}{U'(t_i^0)}=\sum_{l=1}^n\frac{m_l}{t_i^0-z_l}\,.
$$
We have
%$$
%\frac{U'(x)}{U(x)}=\sum_{i=1}^b\frac1{x-t_i^0}+\sum_{l=a}^n\frac{m_l+1}{x-z_l}\,,
%$$
%therefore 
$$
U'(x)=U(x)\left(\sum_{i=1}^b\frac1{x-t_i^0}+\sum_{l=a}^n\frac{m_l+1}{x-z_l}\right)\,,
$$
%and
%$$
%\frac{U''(x)}{U(x)}= \left(\sum_{i=1}^b\frac1{x-t_i^0}+\sum_{l=a}^n\frac{m_l+1}{x-z_l}
%\right)^2
%-\sum_{i=1}^b\frac1{(x-t_i^0)^2}-\sum_{l=a}^n\frac{m_l+1}{(x-z_l)^2}\,,
%$$
%therefore
\begin{eqnarray*}        
U''(x) & = & U(x)   \left(\sum_{i < j}\frac2{(x-t_i^0)(x-t_j^0)}+
\sum_{i=1}^b\sum_{l=a}^n\frac{2(m_l+1)}{(x-t_i^0)(x-z_l)}\right.\\ 
 & + &\left. \sum_{l < j}\frac{2(m_l+1)(m_j+1)}{(x-z_l)(x-z_j)}
+\sum_{l=a}^n\frac{(m_l+1)\,m_l}{(x-z_l)^2}\right)\,.
\end{eqnarray*}        
Thus 
$$
\frac{U''(t_i^0)}{U'(t_i^0)}=\sum_{j\neq i}\frac2{t_i^0-t_j^0}
+\sum_{l=a}^n\frac{2(m_l+1)}{(t_i^0-z_l)}\,,
$$  
hence $(t_1^0,\dots,t_b^0)$ is a solution of the critical point system
of the function  $\Phi_{b, n}(t; z, \tilde m)$. 

\medskip                                                                          
To prove the second statement we 
have to check that $H(x) \, = \, - \, [\, F(x) \, U'' (x) + G (x) \,
U'(x) \, ] \,/ \, U(x)$ is a polynomial.
The requirement that the function $H(x)$ does not have poles at $x
\, = \, t_1^0,\, ...\, , \, t_b^0$ 
is equivalent to the fact that $t^0$ is a critical point of 
$\Phi_{b, n}(t; z, \tilde m)$. 
An easy direct calculation shows that $H(x)$ does not have poles at $x \, =
\, z_a, \, ... \, , \, z_n$.
 \hfill $\triangleleft$

\medskip\noindent 
\begin{lemma}\label{tilde}
Let $m\in \Z^n_{>0},\, \ k\in \Z_{>0}$. Let $\tilde m$ be as in Lemma \ref{ps}.
Assume that $l(m) -2k \geq 0$. Then the pair $\{\tilde m, b\}$ is good.
\hfill $\triangleleft$
\end{lemma}

\medskip\noindent

\begin{lemma}\label{log-min} \ \
Let $m\in \Z^n_{>0}$ and  $ k\in \Z_{>0}$ be such that $k < l(m)+1-k$. 
Let $a\in \Z_{>0}$ be such that $a \leq n$ and $k = (m_a+1)+...+(m_n+1)$.
For pairwise distinct $z_1, ... , z_n$,
consider the differential equation (\ref{Ee}) 
%$$
%F(x)u''(x) + G(x)u'(x) + H(x) u(x) = 0
%$$
where $F(x)$,  $G(x)$ are defined by (\ref{M})
and $H(x)$ is such that the differential equation 
has a solution $U(x) = (x-z_a)^{m_a+1} \cdots (x-z_n)^{m_n+1}$.
Then for  generic $z_1,  ... , z_n$,
the generic solution of this differential equation is multivalued.
\end{lemma}                                                                      

\medskip\noindent
{\bf Proof.}\ \ 
The substitution $u(x) \,= \, U(x)\, v(x)$ turns 
the equation into the differential equation
$$
v''(x) \, + \, \left (\,  \sum_{l=a}^n    {m_l+2\over x - z_l}  
\, - \, \sum_{l=1}^{a-1} {m_l\over x - z_l} \, \right) \,v'(x)\, = \, 0\,.
$$
Its general solution is
$$
v(x) \,=\, \int \,
{ (x-z_1)^{m_1} \, \cdots \, (x-z_{a-1})^{m_{a-1}}
\over
(x-z_a)^{m_a+2}\, \cdots \, (x-z_n)^{m_n+2} }\, \ \ dx\,.
$$
According to our assumptions, \ $m_1 + ... + m_{a-1} \geq
 (m_a + 2) + ... + (m_{n} + 2 ) -1$. In this case the function $v(x)$ is multivalued for generic
$z_1, ... , z_n$. To see this it is enough to notice that the residue of the integrand at
infinity is not zero if $z_a = ... = z_n = 0$
and $z_1 = ... = z_{a-1} = 1$.
 \hfill $\triangleleft$

\medskip\noindent

\begin{thm}\label{log} \ \
Let $m\in \Z^n_{>0}$ and  $ k\in \Z_{>0}$ be such that $k < l(m)+1-k$. 
Let $a\in \Z_{>0}$ be such that $a \leq n$ and $k > (m_a+1)+...+(m_n+1)$.
Set $b= k - (m_a+1)-...-(m_n+1) $. For $s>0$, set
$z^{(s)}=(s\,,\ s^2\,,\ \dots\,, \ s^n)$. Let $t^0$ be any critical point of
the function $\Phi_{b,n}(t, z^{(s)}, \tilde m)$ and let
$E(t^0, z^{(s)}, m, \tilde m)$ be
the associated differential equation. If $s\gg 1$, then the generic solution of
$E(t^0, z^{(s)}, m, \tilde m)$ is a multivalued function.
\end{thm}                                                                      

Theorem \ref{log} implies Theorem \ref{Prop}.

\subsection{Proof of Theorem \ref{log}} \label{PNI}

The pair $\{\tilde m, b\}$ is good by Lemma \ref{tilde}. 
The critical points of the function $\Phi_{b,n}(t, z^{(s)}, \tilde m)$
are labeled by admissible sequences $I=(i_1,...,i_n)$, where $i_1=0$ and
$i_2+...+i_n = b$,  see Sec. \ref{bound}.

Let $t_I(s)=(t_{I,1}(s), ... , t_{I,b}(s))$ be the critical point corresponding
to a sequence $I$. According to the construction, as $s$ tends to infinity for any $j$
there exists the limit of $t_{I,j}(s)/s^n$. This limit is equal to zero if
$j \leq i_1+...+i_{n-1}$,  and the limit is not equal to zero otherwise. Moreover, the limits
of the last $i_n$ coordinates form a critical point of the function
$\Phi_{i_n,\, 2}\, (\, t; \ (0, 1), \ (m_{i_n,1 },  m_{i_n, 2})\, )$ where
$m_{i_n,1}=l(m) - 2k + m_n + 2 + 2 i_n$ and $m_{i_n, 2}= -m_n -2$. 
We denote $(T_1,... , \,T_{i_n})$ the coordinates of that critical point.

Consider the polynomial (\ref{mr}) and make the change of variables $x=s^n y$,
then
\begin{eqnarray*}   
U(s^n y) & = & s^{nk} (y-t^0_1(s)/s^n) \cdots (y-t^0_b(s)/s^n)
\cdot ( y - s^{a-n})^{m_a+1} \cdots ( y - s^{n-n})^{m_n+1} 
\\ 
& = & s^{nk} \,y^{k - i_n - m_n - 1}\, (y-1)^{m_n+1}\,
(y - T_1) \cdots (y - T_{i_n}) + \mathcal O (s^{nk - 1})\,.
\end{eqnarray*}        
Denote 
$$
V(y) \ = \ y^{k - i_n - m_n - 1}\ (y-1)^{m_n + 1}\ (y - T_1) \cdots (y - T_{i_n})\,.
$$

Make the change of variables $x=s^n y$ in 
the differential equation $E( t_I(s), z^{(s)}, m, \tilde m )$,
$$
F(x) u''(x) + G(x) u'(x) + H(x) u(x) = 0\,.
$$
We have 
\begin{eqnarray*}   
F(s^n y) \, &=&\,  s^{n^2} \,(y - s^{1-n}) \cdots (y - s^{n-n}) \,= \,s^{n^2} \,y^{n-1}\, (y-1)\,
\ +  \ \mathcal O (s^{n^2 - 1})\,,
\\
G(s^ny)& =& \ - \ \left( \ \sum_{l=1}^n { m_l\over s^ny - s^l}\ \right) \ F(s^n y)
\\
& = & \ - \, s^{n (n-1)} \,\ \left( \ { l(m) - m_n \over y}\ +
\ {m_{n}\over y - 1} \ \right) \ y^{n-1} (y-1) \ + \ \mathcal O (s^{n (n-1) - 1})\,, 
\\
H(s^ny)\  & = & \- \ {F(s^ny)\ U''(s^n y)\ + \ G(s^n y)\  U'(s^n y)
\over U(s^ny) }\\ 
&= &\ - \ s^{n(n-2)} \  { f(y)\ V''(y)\ + \ g(y) \ V'(y)
\over V(y)} \ +\ \mathcal O (s^{n(n-2) - 1})\,.
\end{eqnarray*}   
Denote
\begin{eqnarray*}   
f(y) & = & \ \ \ \ y^{n-1} (y-1), \\
 g(y) &  =& -  \ \left(\, \ { l(m) - m_n \over y}\, +\,
{m_{n}\over y - 1} \ \right) \ y^{n-1} (y-1)\ , \\
\quad
h(y) & = & -\  \ {f(y) \ V''(y)\  + \ g(y)\  V'(y)
\over V(y)}\,.
\end{eqnarray*}

As  $s\rightarrow\infty\,$, the equation
\begin{eqnarray}\label{s}
F(s^ny)\, u''(s^ny) \,+ \,G(s^ny)\, u'(s^ny)\, +\, H(s^ny)\, u(s^ny) \,= \,0\,
\end{eqnarray}
turns into the equation
\begin{eqnarray}\label{pq}
f(y)\, v''(y)\, +\, g(y)\, v'(y)\, + \,h(y)\, v(y)\, = \, 0\,,
\end{eqnarray} 
and $V(y)$ is its solution. Rewrite  equation  (\ref{pq})
in the form  
$$
v''(y) \, +\, p(y)\, v'(y)\, + \,q(y)\, v(y)\, = \, 0\,.
$$ 
We have
$$
p(y)\, =\, - \, {l(m)- m_n \over y}
\, - \, {m_{n}\over y - 1}\,.
$$

\medskip\noindent
\begin{lemma}\label{q} We have
\begin{eqnarray*}
q(y)& = & -\ \sum_{i < j}\ \frac2{(y-T_i)\ (y-T_j)}-
\sum_{j=1}^{i_n}\ \frac{2\ (k-i_n-m_n-1)}{(y-T_j)\ y}\\
& - & \sum_{j=1}^{i_n}\ \frac{2\ (m_n+1)}{(y-T_j)\ (y-1)}
\ - \ \frac{2\ (k-i_n-m_n-1)\ (m_n+1)}{y \ (y-1)}\\
& - & \ \frac{(k-i_n-m_n-1)\ (k-i_n-m_n-2)}{y^2}
\ -\ \frac{(m_n+1)\ m_n}{(y-1)^2}\\
& + & \left({l(m)- m_n \over y}+{m_n\over y-1}\right)
\ \left(\ \sum_{j=1}^{i_n}\ \frac1{y-T_j}\ +\ \frac{k-i_n-m_n-1}y\ +\
\frac{m_n+1}{y-1}\ \right) \,.   
\end{eqnarray*} 
\end{lemma}      

\medskip\noindent                                                                
{\bf Proof.}\ \ We have
$$
V'(y) \ = \ V(y) \ \left( \ \sum_{j=1}^{i_n}\frac1{y-T_j}+\frac{k-i_n-m_n-1}y+
\frac{m_n+1}{y-1} \ \right)\,,
$$          
and
\begin{eqnarray*}
V''(y) &\  = \ & V(y) \   \left( \ \sum_{i < j}\ \frac2{(y-T_i)\ (y-T_j)}\ +\
\sum_{j=1}^{i_n}\ \frac{2\ ( k-i_n-m_n-1 )}{(y-T_j)\ y}\ \right.\\
& + & \ \sum_{j=1}^{i_n}\ \frac{2\ (m_n+1)}{(y-T_j)\ (y-1)} 
\ + \ \frac{2\ (k-i_n-m_n-1)\ (m_n+1)}{y\ (y-1)}\\
& + & \left.  \frac{(k-i_n-m_n-1)\ (k-i_n-m_n-2)}{y^2}
+ \frac{(m_n+1) \ m_n}{(y-1)^2}\ \right)\,.
\end{eqnarray*}    
Therefore
\begin{eqnarray*}   
q(y)& = & {h(y)\over f(y)} = -\frac{f(y)V''(y)+g(y)V'(y)}{V(y)f(y)}=
-\frac{V''(y)}{V(y)}-p(y)\frac{V'(y)}{V(y)}\\
& = & -\sum_{i < j}\frac2{(y-T_i)(y-T_j)}
-\sum_{j=1}^{i_n}\frac{2(k-i_n-m_n-1)}{(y-T_j)y}\\
& - & \sum_{j=1}^{i_n}\frac{2(m_n+1)}{(y-T_j)(y-1)} 
- \frac{2(k-i_n-m_n-1)(m_n+1)}{y(y-1)}\\
& - & \frac{(k-i_n-m_n-1)(k-i_n-m_n-2)}{y^2}
-\frac{(m_n+1)m_n}{(y-1)^2}\\
& + & \left({l(m)- m_n \over y}+{m_n\over y-1}\right)
\left(\ \sum_{j=1}^{i_n}\ \frac1{y-T_j}\ +\ \frac{k-i_n-m_n-1}y+
\frac{m_n+1}{y-1}\right)\,.             
\end{eqnarray*}

\medskip\noindent 
\begin{lemma}\ \
Equation  (\ref{pq}) is the Fuchsian differential equation
with singular points $0,\, 1,\, \infty$ and
exponents \ $( \,k - i_n - m_n -1,  \ l(m) - k + i_n + 2 \,)\,, \
(\, 0, \,m_n + 1\, ), \ ( \,- k,\, k - l(m) - 1 \,)$, respectively.
\end{lemma}

\medskip\noindent                                                                
{\bf Proof.}\ \ First, we prove that  equation (\ref{pq}) is Fuchsian.
To show this one needs to check that the function $q(y)$ can  be written in the form
$$
\sum\frac {A_i}{y-B_i}+\sum\frac {C_i}{(y-B_i)^2}\,,
$$
where  $B=(0,1,T_1,\dots, T_{i_n})$ and the numbers $A_i$ satisfy the condition
$\sum A_i=0\,$.
This statement clearly follows from the formula for $q(y)$ in Lemma \ref{q}. 
Thus  $p(y)$ and $q(y)$ are of the required form, 
\cite{R}, Ch.~6.41, Theorem 25.     

Now we check that any point $T_j$ is an ordinary point of equation (\ref{pq}).
The formula for $p(y)$ tells that $p(y)$ is holomorphic at $T_j\,$.
To show that $q(y)$  is holomorphic at $T_j\,$, it is enough to verify 
that the limits $q_0$ and $q_1$ of the functions $(y-T_j)^2 q(y)$ and  
$(y-T_j) q(y)$ as $y\rightarrow T_j$ vanish.
                               
All summands in $q(y)$ 
contain $(y-T_j)$ in degree at most $-1\,$. Thus $q_0 = 0\,$. 

We have
$$
q_1\ =\ -\ \frac{m_n+2}{T_j-1}\ +\ \frac{l(m)-2k+m_n+2+2i_n}{T_j}
\ -\ \sum_{i\neq j}\frac2{T_j-T_i}\,.
$$
Thus  $q_1=0\,$, since it is exactly the $j$-th equation 
of the critical point equations for the function
$\Phi_{i_n,2}(\,t;\ (0, 1),\ (m_{i_n,1 }, m_{i_n, 2})\,)\,$. 

The exponents at the singular points $0\,,\ 1\,,\ \infty\,$ are calculated using the indicial
equation and formulas for $p(y), \, q(y)$.
 \hfill $\triangleleft$

\medskip\noindent 
\begin{lemma}\ \
The generic solution of equation  (\ref{pq}) is 
multivalued.
\end{lemma}

\medskip\noindent                                                                
{\bf Proof.}\ \ 
We write the $P$-symbol of equation (\ref{pq}),  see \cite{R}, Ch.~6.45,
$$
v=P\left(\begin{array}{cccc}
         0          &    1     &     \infty    & \      \\
  k - i_n - m_n -1 &    0     &        -k      & \    y \\
  l(m)-k + i_n + 2 & m_n+1    &\   k-l(m)-1    & \
\end{array}\right)\,.
$$
By a meromorphic change of variables,  equation  (\ref{pq})  can be  reduced to the form
\begin{eqnarray}\label{hg}
y (1-y)\ V''(y)\ +\ [\ c - (a+b+1) y\ ]\ V'(y)\ -\ a b\ V( y ) \ =\ 0\,,
\end{eqnarray}
where
$$
a = 2k - l(m) - m_n - i_n - 2\,,
\qquad
b = - i_n - m_n - 1, \qquad 
c = 2k - l(m) - m_n - 2 i_n - 2\,,
$$
see \cite{R}, Ch.~6.46,  \cite{H}, Ch.~2.1.1.

\medskip
Equation (\ref{hg}) is the Gauss hypergeometric equation.
If at least one of the numbers $a\,,\ b\,,\ c-a\,,\ c-b\,$
is an integer, then  formulas for two linearly independent  solutions are 
listed in Ch.~2.1.2 in \cite{H}.  The corresponding table in \cite{H}, pp.~71-73,
consists of  $29$ cases. Moreover,  if generic solutions
are multivalued, then this is stated in the table.

In equation (\ref{hg}), the numbers $a, \ b,\ c$ are
negative integers with $b \geq a \geq c$. This is Case 24 with multivalued 
generic solutions.
\hfill $\triangleleft$

\medskip\noindent                                                                
Let $W(y)$ be a multivalued solution of (\ref{pq}) with  
the initial values at some  point $z_0$ being $W(z_0)=c_1\,,\ W'(z_0)=c_2 \,$
for suitable numbers $c_1, c_2$.
Analytical continuation of this solution 
along some closed curve $\Gamma$ leads to a new value of $W(y)$ at $z_0\,$ 
which is different from the initial value.

Let $X(y,s)$ be the solution of  equation (\ref{s}) with the same initial
values. Then the function $X(y,s)$, restricted to the curve ${\Gamma}$,
tends to the function $W(y)$, restricted to the curve ${\Gamma}$,
as $s\rightarrow\infty\,$. Thus $X(y,s)$ is a multivalued function for  $s \gg 1$.

 This proves  Theorem \ref{log}. \hfill $\triangleleft$

\end{document}